\documentclass[a4paper, reqno, 12pt]{amsart}

\usepackage{geometry}       

\usepackage{float}
\usepackage{bm}
\usepackage{fullpage,xcolor}
\usepackage[mathscr]{euscript}
\usepackage[all]{xy}
\usepackage{epsfig}
\usepackage{amsfonts}
\usepackage{mathptmx}

\usepackage{graphicx}
\usepackage{amssymb} 
\usepackage{amsmath}
\usepackage{amsthm}
\usepackage{mathrsfs}
\usepackage{epstopdf}
\usepackage{url}
\usepackage[msc-links,alphabetic]{amsrefs}
\usepackage{tikz}

\usepackage{color}
\makeatletter

\@addtoreset{equation}{section}
\makeatother 

\theoremstyle{plain}
\newtheorem{theorem}{Theorem}[section]
\newtheorem{corollary}[theorem]{Corollary}

\theoremstyle{definition}
\newtheorem{definition}[theorem]{Definition}

\newtheorem{remark}[theorem]{Remark}

\usepackage{latexsym}
 
\title[The Theory of Doubly Periodic Pseudo Tangles]{The Theory of Doubly Periodic Pseudo Tangles}

\author{Ioannis Diamantis}
\address{Department of Data Analytics and Digitalisation,
Maastricht University, School of Business and Economics,
P.O.Box 616, 6200 MD, Maastricht,
The Netherlands.}
\email{i.diamantis@maastrichtuniversity.nl}

\author{Sofia Lambropoulou}
\address{School of Applied Mathematical and Physical Sciences, National Technical University of Athens, Zografou campus, GR-15780 Athens, Greece.}
\email{sofia@math.ntua.gr}
\urladdr{http://www.math.ntua.gr/~sofia}

\author{Sonia Mahmoudi}
\address{Advanced Institute for Materials Research, Tohoku University, 2-1-1 Katahira, Aoba-ku, Sendai 980-8577, Japan; RIKEN iTHEMS, 2-1 Hirosawa, Wako, Saitama 351-0198, Japan}
\email{sonia.mahmoudi@tohoku.ac.jp}



\date{}

\begin{document}

\setcounter{section}{-1}

\begin{abstract} 
Doubly periodic tangles (DP tangles) are configurations of curves embedded in the thickened plane, invariant under translations in two transversal directions. In this paper we extend the classical theory of DP tangles by introducing the theory of {\it doubly periodic pseudo tangles} (pseudo DP tangles), which incorporate undetermined crossings called {\it precrossings}, inspired by the theory of pseudo knots. Pseudo DP tangles are defined as liftings of spatial pseudo links in the thickened torus, called {\it pseudo motifs}, and are analyzed through diagrammatic methods that account for both local and global isotopies. We emphasize on {\it pseudo cover equivalence}, a concept defining equivalence between finite covers of pseudo motif diagrams. We investigate the notion of equivalence for these structures, leading to an analogue of the Reidemeister theorem for pseudo DP tangles. Furthermore, we address the complexities introduced by pseudo cover equivalence in defining minimal pseudo motif diagrams. This work contributes to the broader understanding of periodic entangled structures and can find applications in diverse fields such as textiles, materials science and crystallography due to their periodic nature.
\end{abstract}

\maketitle


\section{Introduction}\label{sec:0}

A \textit{doubly periodic tangle (DP tangle)} is a set of curves embedded in the thickened plane, invariant under translations in two transversal directions. So, a DP tangle can be defined as the the preimage of a link, called \textit{motif}, embedded in the thickened torus, $T^2\times I$, via a lifting of  $T^2\times I$ under a universal covering map. The topological study of DP tangles can find applications in many fields of science,  such as textiles,  materials science and crystallography  \cites{Grishanov1, GrishanovP1,GrishanovP2,Yaghi,Treacy},  as the mathematical structure of a material or system is often related to its physical properties.  

\begin{figure}[H]
\centerline{\includegraphics[width=5.5in]{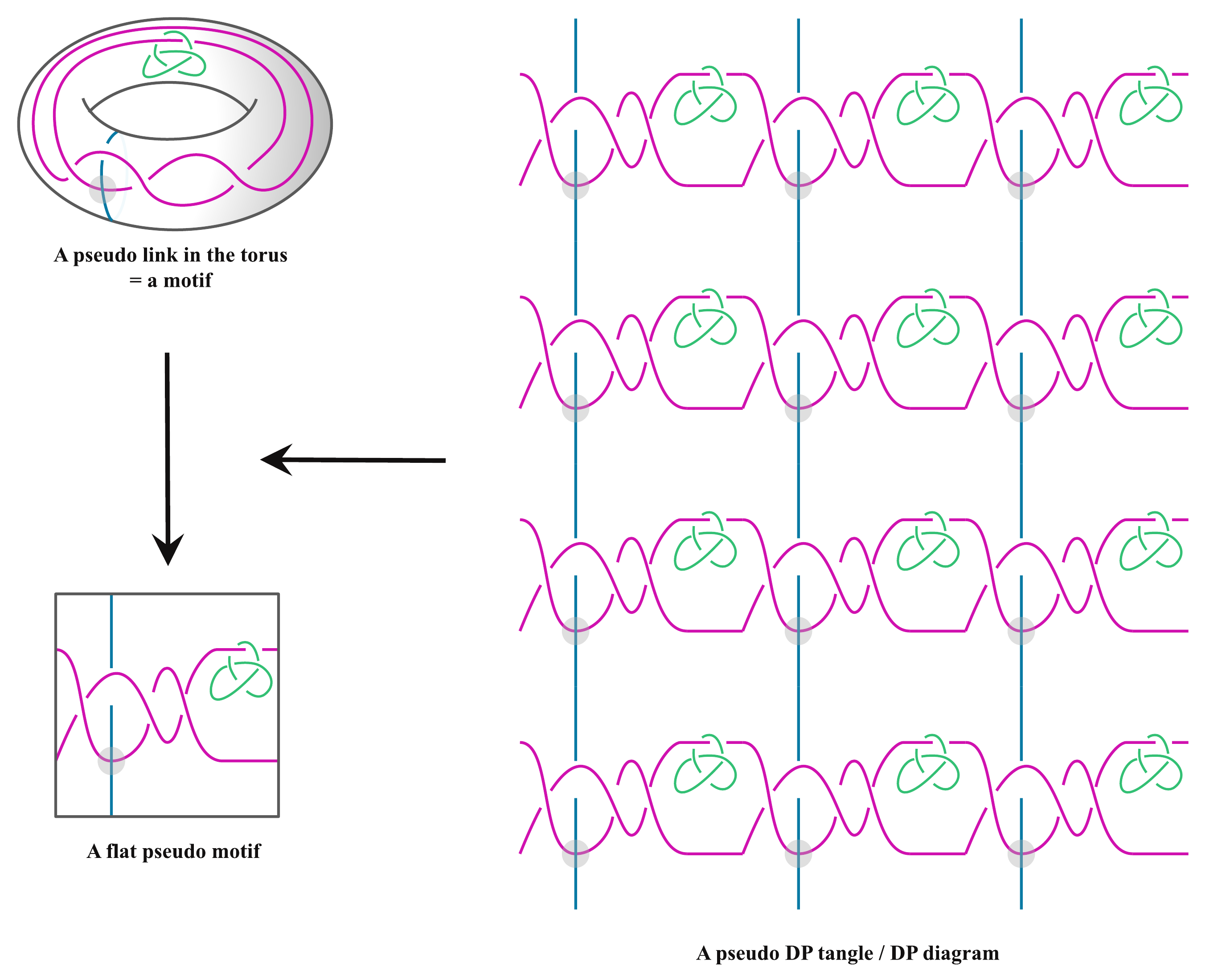}}
\vspace*{8pt} 
\caption{\label{DPtangle} 
A pseudo DP diagram as the lift of a toroidal pseudo link diagram, and a corresponding flat pseudo motif for a choice of an $(l,m)$ pair. }
\end{figure}

As in classical knot theory, the properties of DP tangles (resp. their motifs) are investigated through a two-dimensional diagrammatic theory, where the relative positions of the curves in the thickened plane (resp. the thickened torus) are described by over and under crossing information. However, there are instances where such crossing information is either missing or intentionally suppressed. This is the essence of the theory of \textit{pseudo knots}, introduced by Hanaki in \cite{H}, which is a diagrammatic theory that extends classical knot theory.  The theory of pseudo knots led us to introduce the theory of \textit{doubly periodic pseudo tangles}, abbreviated as {\it pseudo DP tangles}, an extension of classical DP tangles that incorporates undetermined crossings known as \textit{precrossings}. An example is illustrated in Figure~\ref{DPtangle}. 

In this paper we introduce the theory of pseudo DP tangles, we investigate the equivalence relation among them translated into  equivalence of their pseudo motif diagrams,  and we define some topological invariants. Our study of pseudo DP tangles is motivated by potential applications in modeling  complex physical entangled systems by imposing periodic boundary conditions, in which some crossings may be undetectable, for example historical or technical worn textiles and polymer melts. Our work also serves as a basis for future investigations.

Classical DP tangles are said to be {\it DP isotopic} if they are related via an ambient isotopy of the thickened plane that preserves double periodicity.  By definition, the topology of a DP tangle is related to the knot theory of the thickened torus, hence also of a thickened identification parallelogram. DP isotopy gives rise to a diagrammatic theory of {\it equivalence on the motif level}, initiated by Grishanov \cites{Grishanov1, GrishanovP1,GrishanovP2}. This comprises  {\it local equivalence}, including  isotopies in the thickened torus, {\it global equivalence}, deriving from orientation preserving affine transformations of the thickened plane (translations, rotations, stretchings, contractions, shearings), and from the equivalence between different finite covers/quotients of the same motif, called \textit{cover equivalence} (or \textit{scale equivalence} \cite{Sonia}). Note that translations correspond to shifts of the longitude-meridian pair of the torus, stretchings and contractions correspond to inflations and deflations of the torus, and shearings and rotations are related to Dehn twists. For a detailed exposition of the DP tangle equivalence  on the level of motifs, we refer the reader to \cite{DLM1} and to references therein.

The theory of pseudo DP tangles is, by definition, built from a diagrammatic theory. In order to study their equivalence relation,  which is the main objective of the paper, we first consider the {\it local pseudo isotopy moves} comprising surface isotopies, the classical Reidemeister moves as well as the pseudo Reidemeister moves involving precrossings, see Figure~\ref{reid}. We further define a pseudo DP tangle as the {\it lift} of a pseudo DP diagram in the thickened plane, using the notion of the lift of a toroidal pseudo link diagram in the thickened torus, as established in \cite{DLM4}, whereby every classical crossing is embedded in a sufficiently small 3-ball, while the precrossings are supported by sufficiently small rigid discs embedded in three-space (see Section~\ref{sec:lift}). Using the lift, we then resource to the theory of \textit{mixed links} (\cite{LR1}), which are links in three-space containing a fixed sublink, since the thickened torus can be described as the complement in three-space of the Hopf link, denoted by ${\rm H}$. So, we extract a complete set of local moves between diagrams of pseudo DP tangles and their motifs. View Figures~\ref{mpr}, ~\ref{mr3} and ~\ref{fig:iso}. 

Finally, we have the notion of {\it pseudo cover equivalence}, which defines an equivalence relation between two finite covers or finite quotients of the same pseudo motif. There is a crucial subtlety to take into consideration in the cover equivalence of pseudo DP tangles. Given a pseudo motif diagram $d$, each precrossing bears the missing information of over/under for each arc, which is unknown to us. This missing information will be inherited by every corresponding precrossing generated by a finite cover or the universal cover of $d$.  This observation leads to the concept of \textit{homologous precrossings}, that is, precrossings lying  in the preimage of the same precrossing  under a covering map. The concept of homologous precrossings makes the notion of pseudo cover equivalence well-defined. Indeed, the converse operation, namely, taking a quotient of $d$, is not well-defined unless its precrossings are labeled as homologous. The above culminate to the following theorem, which is the main result of the paper.

\begin{theorem}[Pseudo DP tangle equivalence] \label{main}
Two pseudo DP tangles are pseudo DP equivalent if and only if any two pseudo (flat) motif diagrams of their are related by a finite sequence of local pseudo motif isotopy moves, longitude-meridian shifts, torus inflations and deflations, Dehn twists, and pseudo cover equivalence.
\end{theorem}

We continue by extending the concepts of resolution set and weighted resolution set for pseudo DP tangles. The \textit{resolution set} of a pseudo DP diagram is defined as the set of all classical DP tangles that can be obtained from a given pseudo DP diagram by resolving homologous precrossings into the same type of crossing. The \textit{weighted resolution set} of a pseudo DP diagram is a generalization of the resolution set that takes also into account the probability that a resolution appears in the resolution set. These two sets are shown to be topological invariants of pseudo DP tangles.

\small
The paper is organized as follows. In \S~\ref{sec:prel} we recall the basics of the theory of  pseudo knots and links in the plane and in the torus. In \S~\ref{sec:psDP} we recall the basic concepts of classical DP tangles and introduce the general framework for pseudo DP tangle diagrams and pseudo DP tangles in the thickened plane. In  \S~\ref{sec:equivalence} we establish the pseudo DP tangle equivalence  on the level of pseudo motifs. Finally, in \S~\ref{sec:resolution}, we define the notions of resolution set and weighted resolution set for pseudo DP tangles.

\section{Planar and toroidal pseudo links}\label{sec:prel}

In this section we recall the basics of the theory of  pseudo knots and links in the plane and in the torus from \cite{DLM4}. We first recall that a  knot is an embedding of a  circle into the three-space, while a  link is an embedding of a finite collection of circles. 

\subsection{Definitions and equivalences}\label{sec:prel-def}

The classical by now theory of planar pseudo knots was introduced by Hanaki in \cite{H} and the mathematical background of (planar) pseudo knot theory was established in \cite{HJMR}. We shall say `link' for both knots and links. Throughout we shall say (pseudo) knots or (pseudo) links referring to both (pseudo) knots and links.

\begin{figure}[H]
\begin{center}
\includegraphics[width=1.2in]{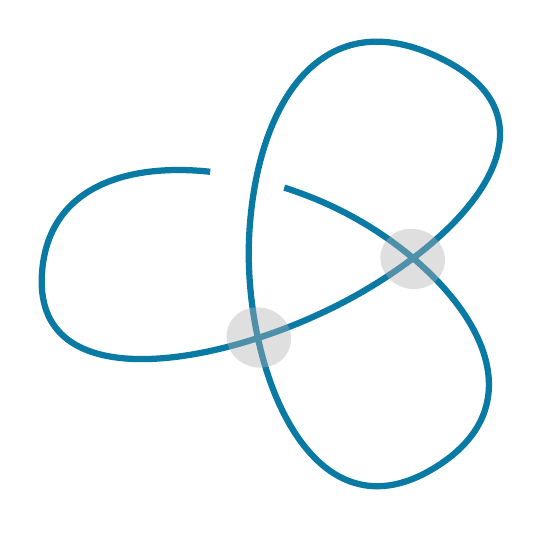}
\end{center}
\caption{A pseudo trefoil.}
\label{pk1}
\end{figure}

A {\it pseudo knot diagram} in an orientable surface $\Sigma$ consists of a regular knot diagram (which is a knot projection with finitely many double points) where some crossing information may be missing, that is, it is unknown to us which arc passes over and which arc  passes under the other. These undetermined crossings are called {\it precrossings} or {\it pseudo crossings} and are depicted as self-intersections of the arcs enclosed in a disc, illustrated in gray  (for an example see Figure~\ref{pk1}). In the present work the surface $\Sigma$ is the plane or the torus.

\begin{figure}[H]
\begin{center}
\includegraphics[width=5.8in]{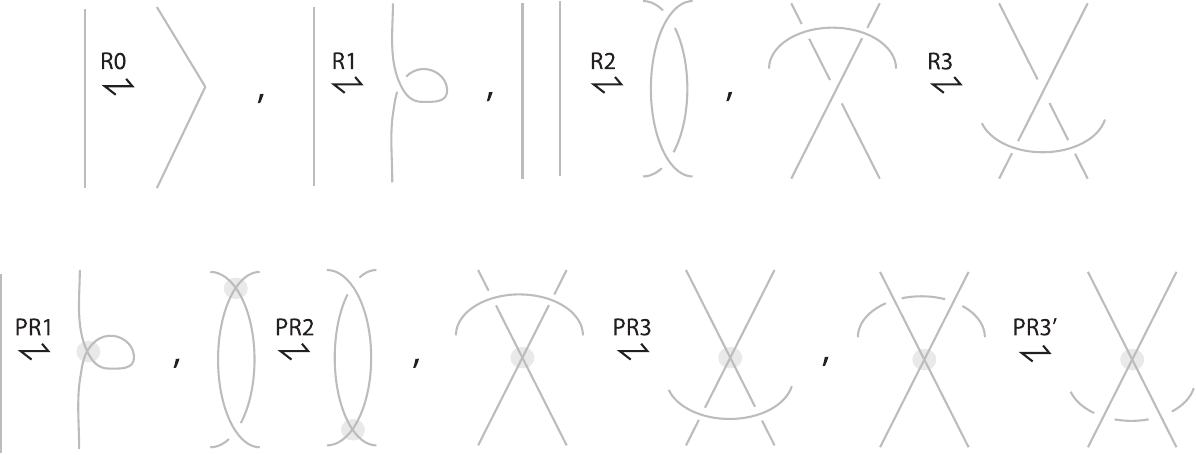}
\end{center}
\caption{Equivalence moves for pseudo knots.}
\label{reid}
\end{figure}

A {\it pseudo knot} in the surface $\Sigma$ is defined as the equivalence class of all pseudo knot diagrams in $\Sigma$ under surface isotopy and all versions of the well-known \textit{Reidemeister moves} for classical knots and their versions containing pseudo crossings, the \textit{pseudo Reidemeister moves}, as exemplified in Figure~\ref{reid}. 

\begin{figure}[H]
\begin{center}
\includegraphics[width=1.7in]{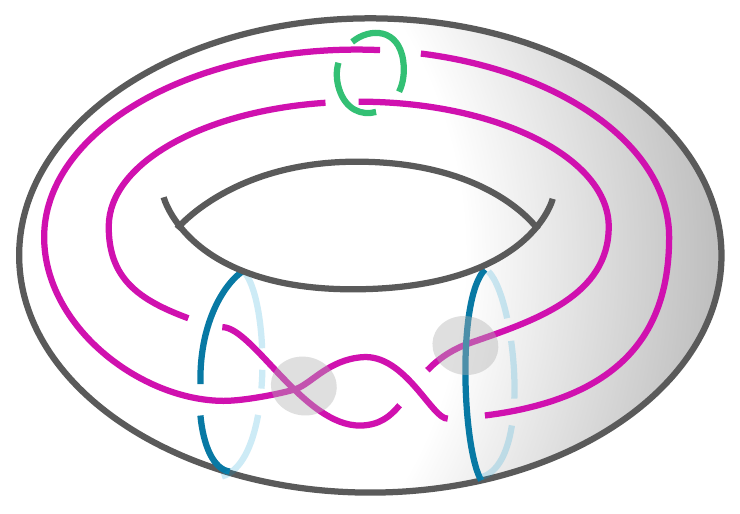}
\end{center}
\caption{A toroidal pseudo link.}
\label{fig:pkannulusttorus}
\end{figure}

In \cite{DLM4} we extended the theory of planar pseudo links to pseudo links in the annulus (cf. also \cite{D}) and in the torus. In the present work we shall only consider toroidal pseudo links. View Figure~\ref{fig:pkannulusttorus} for an example.

\subsection{Lifts of planar and toroidal pseudo links}\label{sec:lift}

\begin{figure}[ht]
\begin{center}
\includegraphics[width=5.5in]{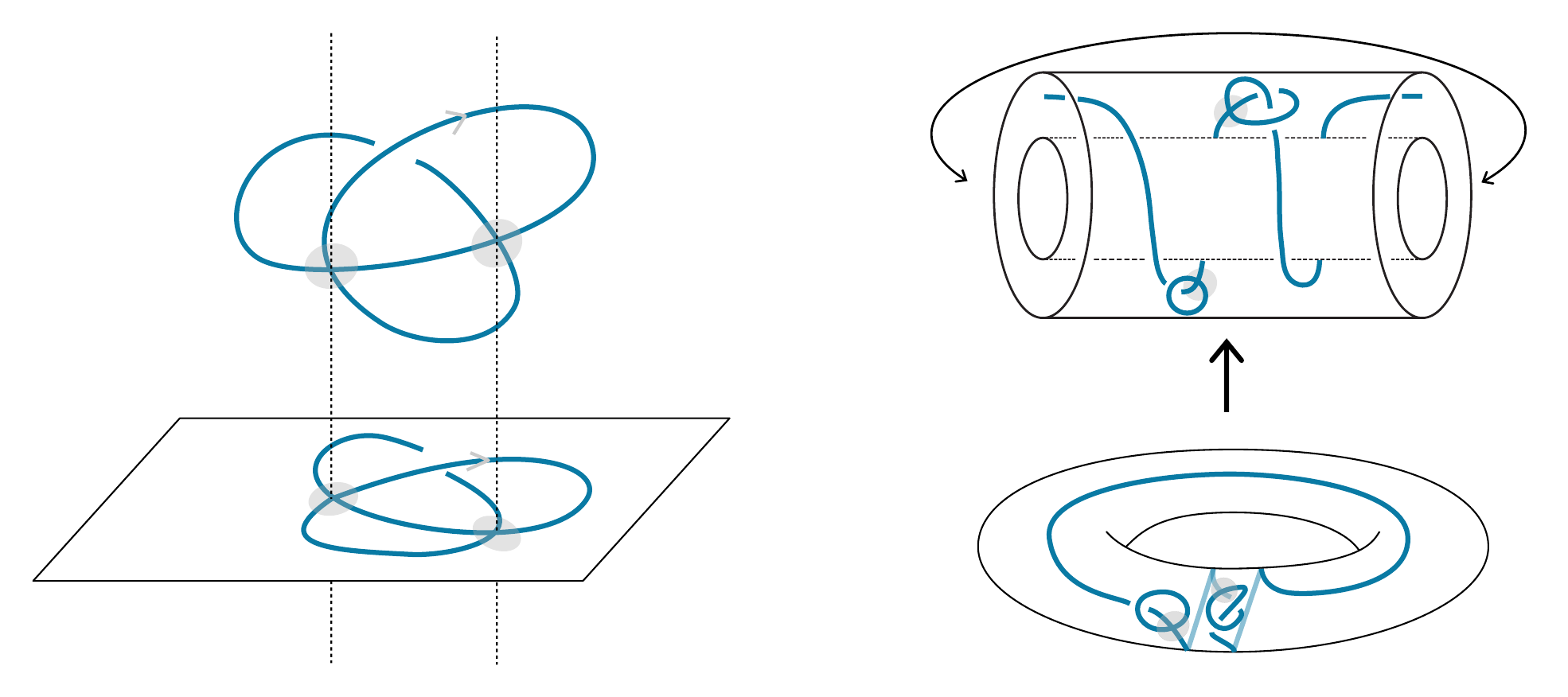}
\end{center}
\caption{Lifting a planar and a toroidal pseudo link.}
\label{fig:alllifts}
\end{figure}

In \cite{DLM4} we first defined the {\it lift} of a planar pseudo link as follows: every classical crossing is embedded in a sufficiently small 3-ball so that the over arc is embedded in its upper boundary and the under arc is embedded in its lower boundary, while the precrossings are supported by sufficiently small rigid discs, which are embedded in three-space.  By lifting the precrossings in rigid discs  we preserve the pseudo link's essential structure while preserving the missing informations of its precrossings.  The simple arcs connecting crossings can be also replaced by isotopic ones in three-space. The resulting lift, is called  {\it spatial pseudo link} and it is a collection of closed curves in three-space consisting in embedded discs from which emanate embedded arcs. For an example view Figure~\ref{fig:alllifts}. The natural equivalence between spatial pseudo links is isotopy. Two (oriented) spatial pseudo links are said to be {\it isotopic} if they  are related  through ambient space orientation preserving homeomorphisms  respecting arc and disc isotopies. In this context, two  (oriented) spatial pseudo links are  isotopic if and only if any two corresponding planar pseudo link diagrams of theirs are (oriented) pseudo Reidemeister equivalent. 

In analogy,  we defined in \cite{DLM4} the lift of a toroidal  pseudo link in the interior of the thickened torus $T^2 \times I$. In Figure~\ref{fig:alllifts} we illustrate the lift of  a toroidal pseudo link in the thickened torus, viewed  as the identification space of a thickened cylinder, by identifying the two annular boundary components. The isotopy of a lift in  $T^2 \times I$ is defined in analogy to the isotopy of a spatial lift and is confined in  the thickened torus. Further, the  isotopy equivalence corresponds bijectively to the toroidal pseudo Reidemeister equivalence. For details see \cite{DLM4}.

\subsection{Mixed pseudo links}\label{sec:mixed}

\begin{figure}[H]
\begin{center}
\includegraphics[width=5.5in]{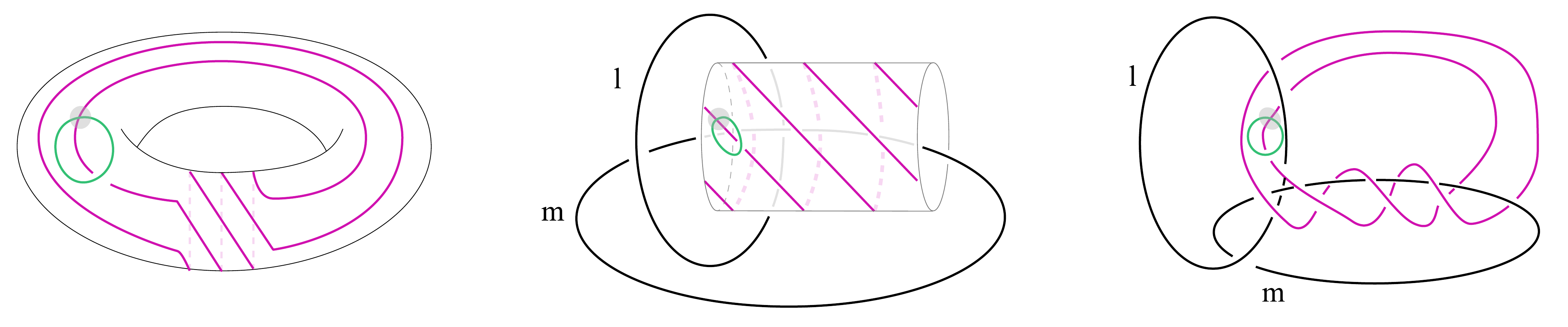}
\end{center}
\caption{A toroidal pseudo link and its corresponding ${\rm H}$-mixed pseudo link.}
\label{mixed}
\end{figure}

The notion of lift for toroidal pseudo knots enables us to establish further connections of this theory to the theory of \textit{mixed links} \cite{LR1}.  It turns out that toroidal pseudo knots are connected to ${\rm H}$-mixed pseudo links in $S^{3}$, which are planar pseudo links that contain a point-wise fixed Hopf link as a sublink, since its complement is a thickened torus. Note that the two fixed components of the Hopf link represent a meridian and a longitude of the torus, respectively. View Figure~\ref{mixed} for an illustration.

\begin{figure}[H]
\begin{center}
\includegraphics[width=5.1in]{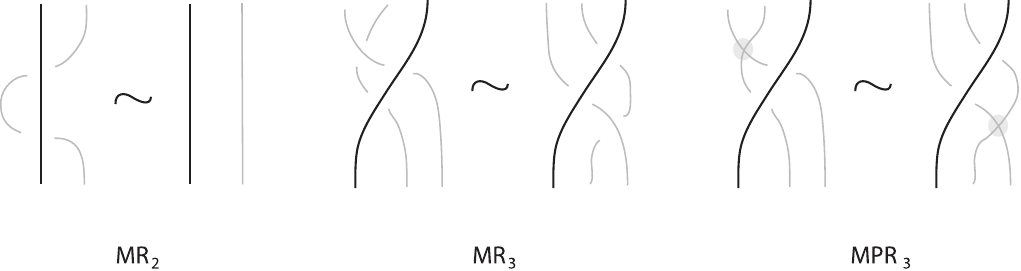}
\end{center}
\caption{Mixed pseudo Reidemeister moves.}
\label{mpr}
\end{figure}

\begin{figure}[H]
\begin{center}
\includegraphics[width=2.3in]{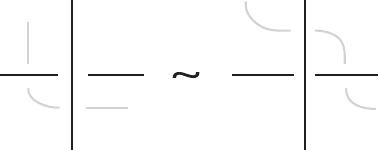}
\end{center}
\caption{An HR3 mixed Reidemeister move.}
\label{mr3}
\end{figure}

As we showed in \cite{DLM4}, isotopy classes of  pseudo links in the thickened torus are in bijection with isotopy classes of ${\rm H}$-mixed pseudo links in $S^{3}$, via isotopies that keep  ${\rm H}$ point-wise fixed. Based on the above, we further translated the spatial mixed pseudo link isotopies to the  diagrammatic level through discrete Reidemeister type moves for  ${\rm H}$-mixed pseudo link diagrams in the plane. These moves preserve the oriented ${\rm H}$ point-wise and include planar isotopies, the local moves for pseudo links, as exemplified in Figure~\ref{reid}  and the mixed pseudo Reidemeister moves exemplified in Figures~\ref{mpr} and~Figure~\ref{mr3}.

\begin{theorem}{(Thm~3.10 \cite{DLM4})}
Two pseudo links in $T^2 \times I$ are isotopic if and only if any corresponding mixed pseudo link diagrams of theirs differ by a finite sequence of planar isotopies and 
\smallbreak
\begin{itemize}
\item[i.] the classical Reidemeister moves R1, R2, R3 and the pseudo Reidemeister moves PR1, PR2, PR3, PR3$^{\prime}$ for the moving parts (view Figure~\ref{reid}),
\smallbreak
\item[ii.] the mixed Reidemeister moves MR2, MR3, MPR3, HR3 that involve the fixed and the moving part of the mixed pseudo diagrams, illustrated in Figures~\ref{mpr} and~\ref{mr3}.
\end{itemize}
\end{theorem}

\begin{remark}\rm
As observed in \cite{DLM4}, by omitting the PR1-moves (recall Figure~\ref{reid}) and by considering the precrossings in Figures~\ref{reid} and~\ref{mpr} as singular crossings, we obtain the analogous  Reidemeister-type  theorem for singular links in $T^2 \times I$.
\end{remark}

\subsection{Resolution and weighted resolution sets}\label{sec:were}

In this subsection we recall an invariant of a pseudo link, known as the {\it  resolution set} \cite{HJMR}. 

\begin{figure}[ht]
    \centering
    \includegraphics[width=0.6\linewidth]{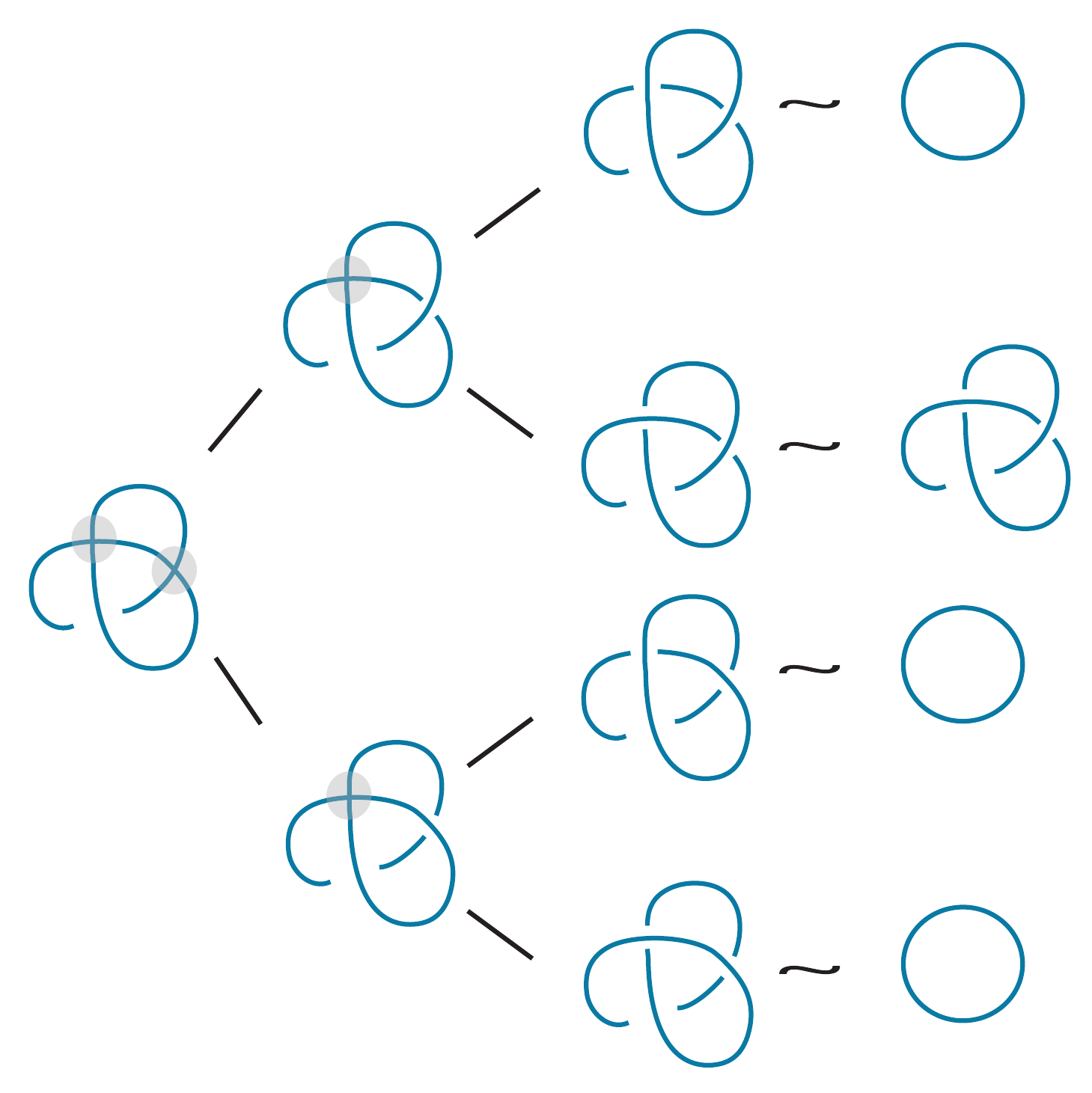}
    \caption{The resolution set of a pseudo knot.}
    \label{fig:weretr}
\end{figure}

A {\it resolution} of a planar pseudo link diagram $K$ refers to a specific assignment of crossings (either over or under) for each precrossing in $K$, resulting in a classical link. The set of all possible resolutions is the {\it resolution set} of $K$. For an illustration see Figure~\ref{fig:weretr}. It is straightforward to verify that the resolution set remains invariant under pseudo link equivalence moves \cite{HJMR}, and hence,  it is an invariant of the pseudo link represented by the diagram $K$. Furthermore, applying any classical invariant to the elements of the resolution set also induces an invariant for the planar pseudo link. The {\it weighted resolution set} or {\it WeRe set} of a planar pseudo link diagram \( K \) is a collection of ordered pairs \((K_i, p_{K_i})\), where \( K_i \) represents a resolution of \( K \), and \( p_{K_i} \) denotes the probability of obtaining from \( K \) the equivalence class of \( K_i \)  through a random assignment of crossing types, with equal likelihood for positive and negative crossings. 

In analogy, we have from \cite{DLM4} the extension of the concept of the resolution set for toroidal pseudo links. 

\begin{definition}\label{wrstormixed}\rm 
A {\it resolution} of a toroidal pseudo link diagram \( K \) is a toroidal  link diagram obtained through a specific assignment of crossing types (positive or negative) for every precrossing in \( K \). The resolution lifts to a link in the thickened torus $T \times I$. The set of all possible resolutions of \( K \) is the {\it toroidal resolution set} of $K$. 
\end{definition}

As in the case of planar pseudo links, in \cite{DLM4}  we showed the following:

\begin{theorem}\label{th:toroidalWeRe}
The toroidal resolution set is an invariant of toroidal pseudo links. 
\end{theorem}

\begin{corollary}
Any invariant of links in the thickened torus applied on the elements of the toroidal resolution set induces an invariant set of the toroidal pseudo link.  
\end{corollary}

\section{The theory of doubly periodic pseudo tangles}\label{sec:psDP}

Inspired by our earlier works, in this section we introduce the theory of {\it pseudo DP tangles}, extending the theory of  DP tangles \cites{Grishanov1,GrishanovP1, GrishanovP2, Sonia, DLM1}. A {\it DP tangle} is the lift to the Euclidean thickened plane  of a link embedded in the thickened torus under the universal covering map. 

In this section we first define the notion of a pseudo DP tangle diagram, since the theory of pseudo knots is a diagrammatic theory, and the notion of generating pseudo motif. A crucial concept in the theory is that of homologous precrossings in a finite cover or the universal cover of a pseudo motif, generated by the same precrossing in the pseudo motif. 
We then define the notion of a pseudo DP tangle, using the lifting technique for precrossings described in the previous section. 
We continue with introducing the notion  of a supporting lattice and we conclude the section with  the notions of minimal lattice and minimal generating pseudo motif. In the next section we discuss the equivalence of pseudo DP tangles.

\subsection{Pseudo DP tangle diagrams}\label{subs:psDPdiagram}

Let $\mathbb{E}^2$ denote the Euclidean plane and let $B=(u,v)$ be a basis of $\mathbb{E}^2$. Let further $\rho$ be the universal covering map $\rho: \mathbb{E}^2 \rightarrow{} T^2$ that assigns a longitude $l$ of $T^2$ to $u$ and a meridian $m$ of $T^2$ to $v$. We then define the following. 

\begin{definition}\label{def:pseudo-DP-tanglediagram}
Let $d$ be a pseudo link diagram in the torus $T^2$. Then, the lifting of $d$ under the covering map $\rho$ is a {\it doubly periodic pseudo  diagram}  (or {\it pseudo DP diagram}) $d_{\infty}$, and  $d$ is  a {\it (generating) pseudo motif} for $d_{\infty}$. One may also view  $d$  as a pseudo tangle diagram in an identification parallelogram (or flat torus), obtained from $T^2$ by cutting along the longitude-meridian pair $(l,m)$,  and this shall be called a {\it  (generating) flat pseudo motif} for $d_{\infty}$. 
\end{definition}

An example of a (flat) pseudo  motif and its corresponding pseudo DP diagram is illustrated in Figure~\ref{DPtangle}. We shall now introduce an important notion in the theory of pseudo DP tangles. 

\begin{definition}\label{def:homologous}
Two precrossings $x_i, x_j$ in a pseudo DP tangle diagram $d_{\infty}$ are said to be  {\it homologous} if they both lie in the preimage of the same precrossing  $x$ of a (flat) pseudo motif $d$  under $\rho$. So, the precrossings $x_i, x_j$ in $d_{\infty}$ bear the same missing information as $x$.  The same definition applies for homologous precrossings $x_i, x_j$ of any finite cover of $d$. 
\end{definition} 

\begin{remark} \label{rem:partition} 
It follows from the definition that the set of all precrossings of $d_{\infty}$  is partitioned into equivalence classes of homologous precrossings.
\end{remark} 

\begin{figure}[ht]
\includegraphics[width=2.7in]{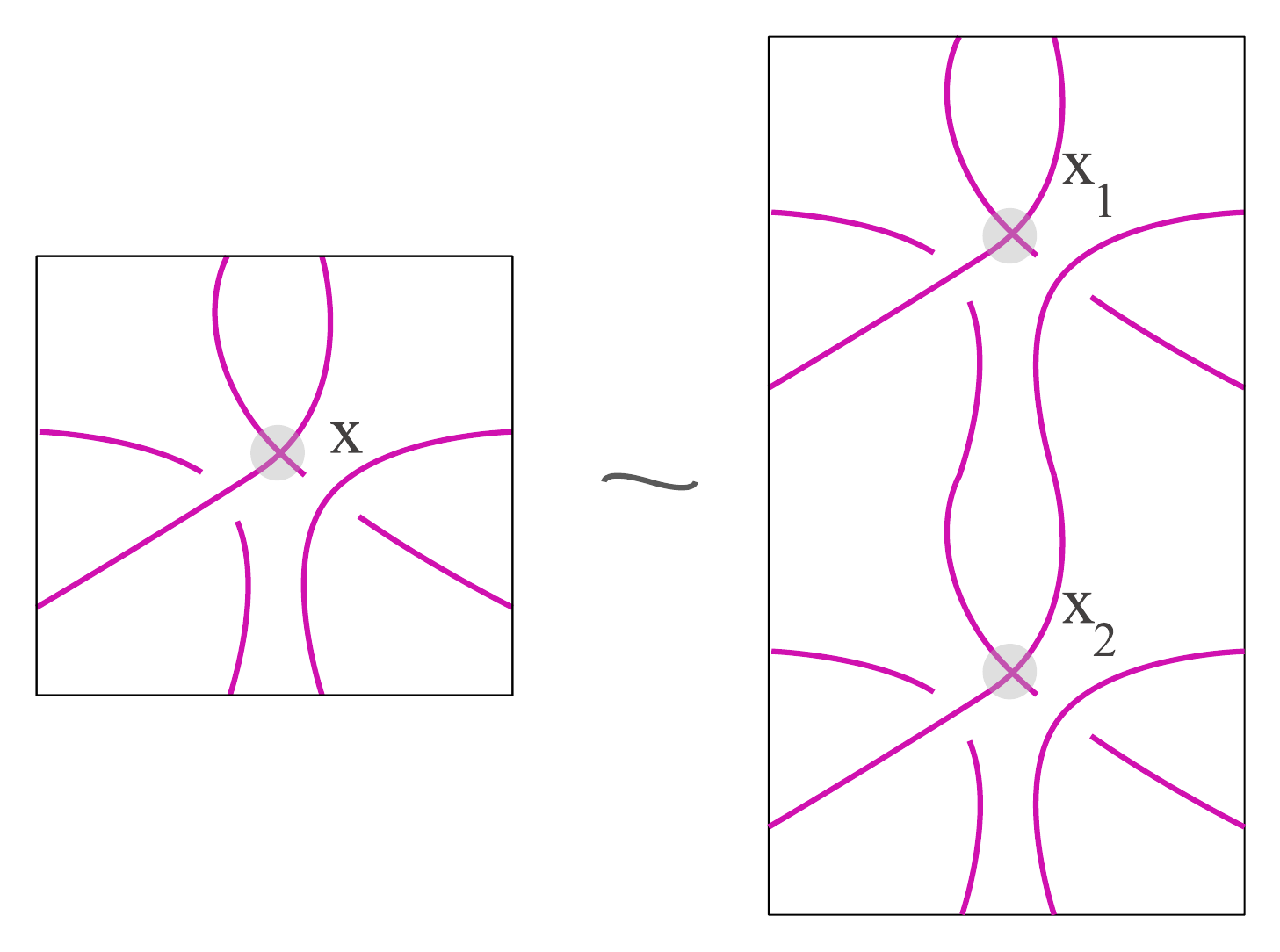}
\caption{\label{Tknot-Tlink1} A  flat pseudo motif on the left and its double cover on the right with two homologous precrossings.}
\end{figure}

\begin{remark} \label{rem:homologous} 
It follows from Definitions~\ref{def:pseudo-DP-tanglediagram} and~\ref{def:homologous}, that any finite cover of $d$ also defines a  (flat) pseudo motif of $d_{\infty}$. View Figure~\ref{Tknot-Tlink1} for an example, where the precrossings $x_1, x_2$ lie in the preimage of $x$, so are homologous. If, however, the two precrossings $x_1, x_2$ in the  right-hand-side illustration of the figure were not marked as homologous, then the right-hand-side motif is not necessarily a motif for the same pseudo DP tangle diagram $d_{\infty}$ as is $d$. 
\end{remark} 

\subsection{Pseudo DP tangles}\label{subs:psDP}

\begin{figure}[ht]
\centerline{\includegraphics[width=3.5in]{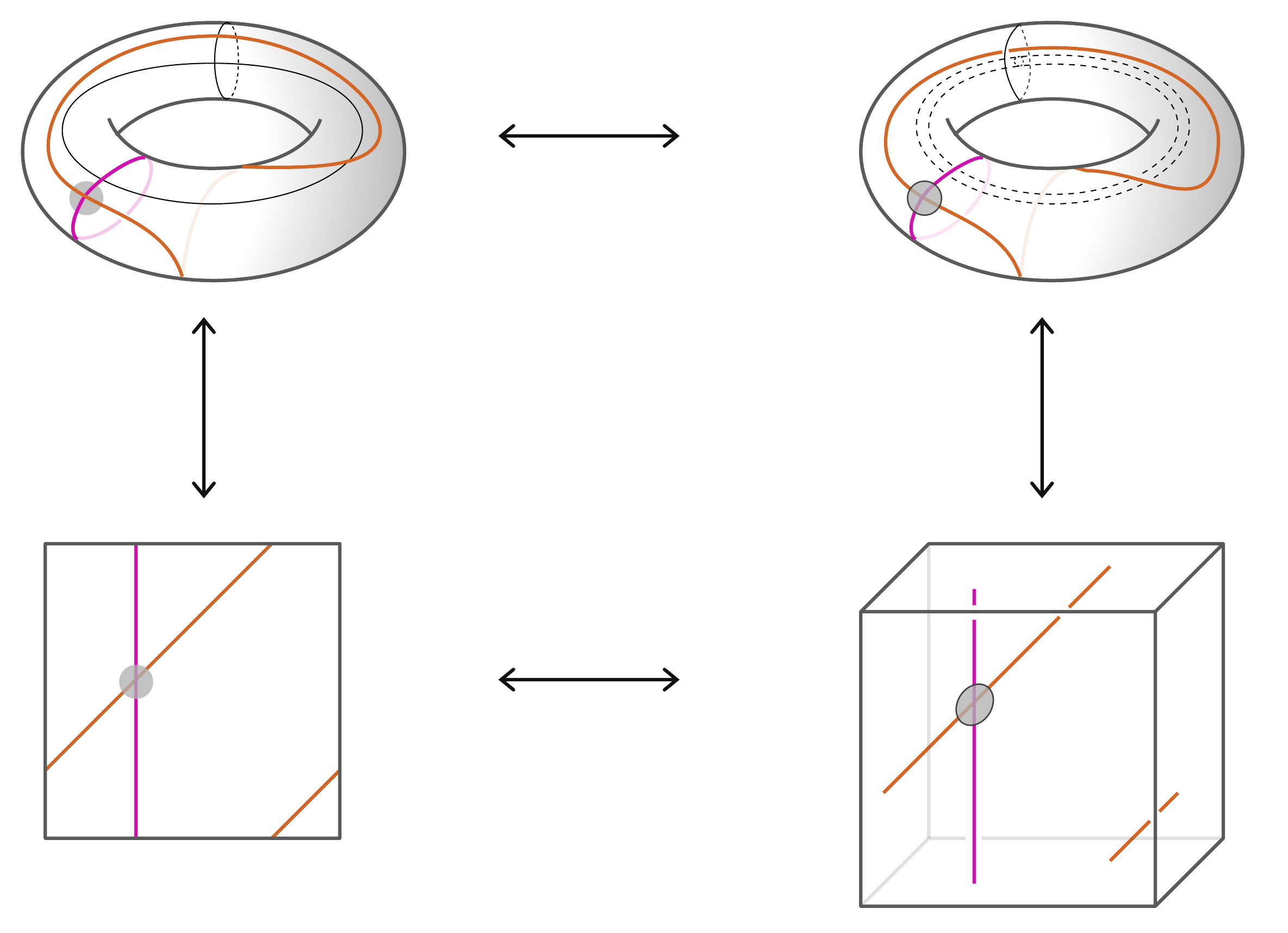}}
\caption{\label{DPlift} 
A toroidal pseudo link diagram  (top left), a corresponding flat motif in a parallelogram (bottom left) and their corresponding spatial lifts in the thickened torus (top right) and in the thickened parallelogram (bottom right), with the embedded precrossing discs.}
\end{figure}

In Section~\ref{sec:lift} we explained that a toroidal pseudo link diagram $d$ has a lift to a toroidal pseudo link $\tau$ in the thickened torus $ T^2 \times I$, where precrossings are lifted as embedded precrossing discs respecting their missing information. Note further that the universal covering map  $\rho$ extends trivially to a universal covering map of the thickened plane to the thickened torus $\tilde{\rho}$: $ \mathbb{E}^2 \times I \rightarrow{} T^2 \times I$. Similarly, the spatial lift of the corresponding to $d$ flat pseudo motif diagram as a pseudo tangle in a thickened parallelogram, subjected to periodic boundary conditions, is well-defined by cutting the thickened torus, bearing the toroidal pseudo link $\tau$,  along the longitudinal ribbon $l \times I$ and the meridional ribbon $m \times I$, as illustrated in Figure~\ref{DPlift}.  Hence, the following notions are well-defined. 
  
\begin{definition}\label{def:pseudo-DP-tangle}
A {\it doubly periodic pseudo tangle} (or {\it pseudo DP tangle}) is defined as the lifting $\tau_{\infty}$  of a toroidal pseudo link $\tau$ in the thickened torus $T^2 \times I$ under the covering map $\tilde{\rho}$, and  $\tau$ is a {\it (generating) pseudo motif} for $\tau_{\infty}$. Further, viewing  $\tau$  as a pseudo tangle in an identification thickened parallelogram (or flat thickened torus), we obtain a {\it  (generating) flat pseudo motif} for $\tau_{\infty}$. 
\end{definition} 

It follows from Definitions~\ref{def:pseudo-DP-tangle} and~\ref{def:pseudo-DP-tangle} that a pseudo DP tangle $\tau_{\infty}$ can be equivalently viewed as the spatial lift of a pseudo DP tangle diagram $d_{\infty}$, where the same missing information is retained on homologous precrossings.
 
\subsection{Periodic lattices and minimal motifs of pseudo DP tangles}\label{sec:lattice}

The set of points $\Lambda (u,v) = \{xu + yv\, |\, x,y \in \mathbb{Z}\} \approx \mathbb{Z}^2$ generated by the basis $B$ of $\mathbb{E}^2$, defines a {\it periodic lattice} for the pseudo DP tangle $\tau_{\infty}$. In particular, it is well-known (see {\cite{DLM1}, \cite{Grishanov1}) that the same periodic lattice $\Lambda = \Lambda (u,v) = \Lambda' (u',v')$ can be generated by two different bases $B=(u,v)$ and $B'=(u',v')$ if and only if for $x_1, x_2, x_3, x_4 \in \mathbb{Z}$,
$$\begin{pmatrix}
u' \\
v' 
\end{pmatrix}
=
\begin{pmatrix}
x_1 & x_2 \\
x_3 & x_4 
\end{pmatrix}
\cdot
\begin{pmatrix}
u \\
v 
\end{pmatrix}
{\it  , \quad where }
\mid{x_1 x_4 - x_2 x_3}\mid = \pm 1.$$

\begin{figure}[H]
\centerline{\includegraphics[width=3.7in]{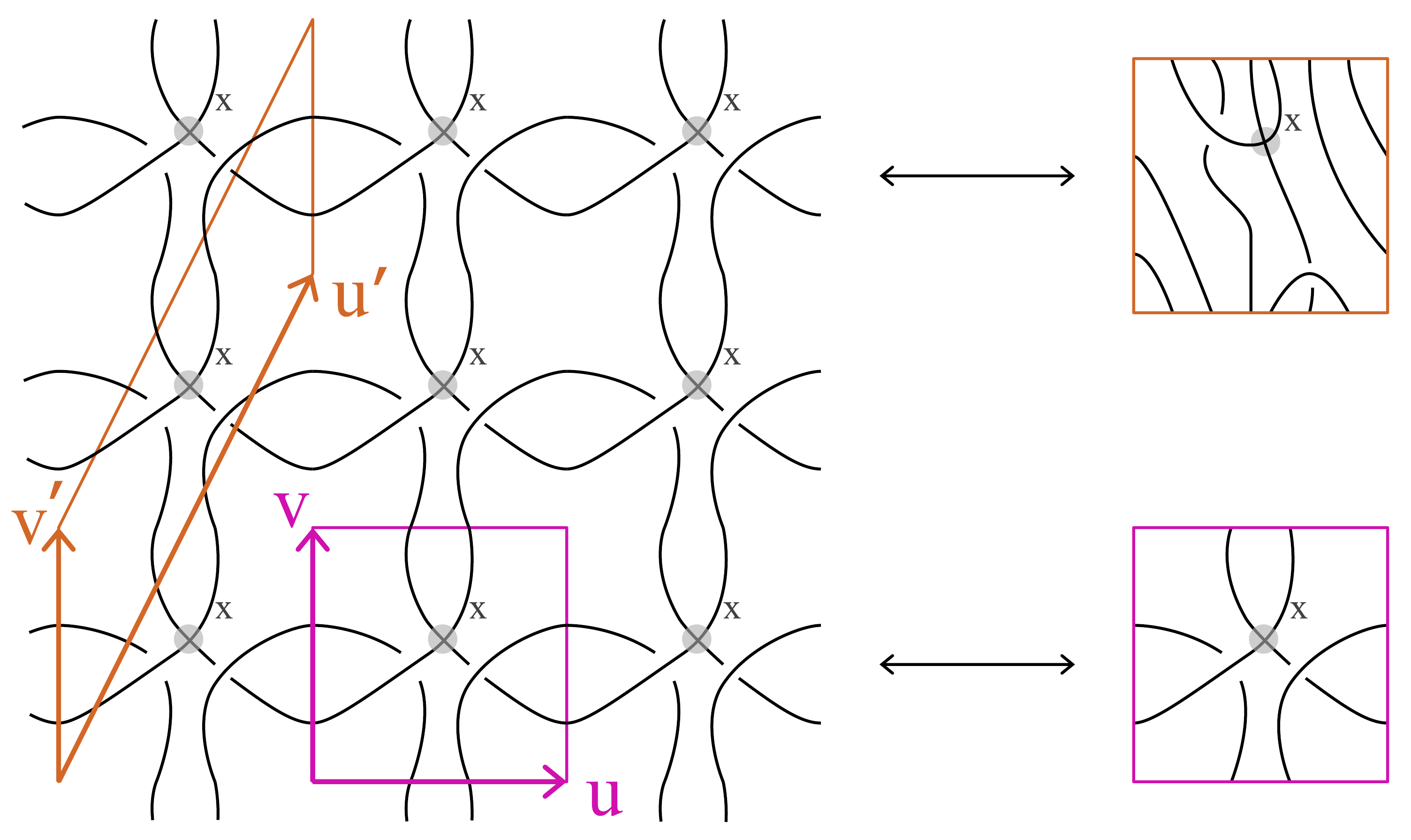}}
\vspace*{5pt}
\caption{\label{shearing} 
A fixed lattice  $\Lambda $ for a pseudo DP tangle, with two different bases of $\mathbb{E}^2$ (orange and purple) and their flat motifs.}
\end{figure}

Let $d$ be a (flat) pseudo  motif  of a pseudo DP tangle diagram $d_{\infty}$ associated with a supporting lattice $\Lambda$.  If  $d$ can be subdivided into a number of identical pieces, according to  Remark~\ref{rem:homologous}, then each one  also defines a  (flat) pseudo motif for $d_{\infty}$. This leads to the definitions of a minimal lattice and minimal motif associated with a pseudo DP tangle.

\begin{definition}\label{def:minimal lattice}
A {\it minimal (flat) motif (diagram) of} $\tau_{\infty}$, resp. $d_{\infty}$, denoted $\tau_{min}$  resp. $d_{min}$, is a motif which is not a finite cover of another (flat) motif (diagram) of $\tau_{\infty}$, resp. $d_{\infty}$. The periodic lattice associated to a minimal (flat) motif diagram is called a {\it minimal lattice}, denoted $\Lambda_{min}$, and we have $d_{min} = d_{\infty} / \Lambda_{min}$.
\end{definition} 

Most topological invariants for DP tangles are based on the assumption of a minimal motif (cf.\cites{Grishanov1,GrishanovP2}).  We note that area-preserving transformations of $\mathbb{E}^2$ retain  lattice minimality.  However, it can be very tricky to find a minimal motif for a DP tangle, and the notion of cover equivalence plays a crucial role in such a search (cf. discussion and example in \cite{DLM1}). In the presence of precrossings this task can be even more subtle if  precrossings are not identified as homologous in a given pseudo motif. In Subsection~\ref{sec:scaleequivalence} we expand on this discussion and present instructive examples in Figures~\ref{tortref} and~\ref{tortref2}. 

Figure~\ref{minimal-pseudo} illustrates an example of a pseudo DP diagram with homologous precrossings denoted by $x$. The two pink flat motifs are expected to be minimal and are related by a Dehn twist. The blue flat pseudo motif is a double cover of the pink one depicted on the bottom right, while the green flat pseudo motif is a double cover of the blue one. 

\begin{figure}[ht]
\includegraphics[width=5.5in]{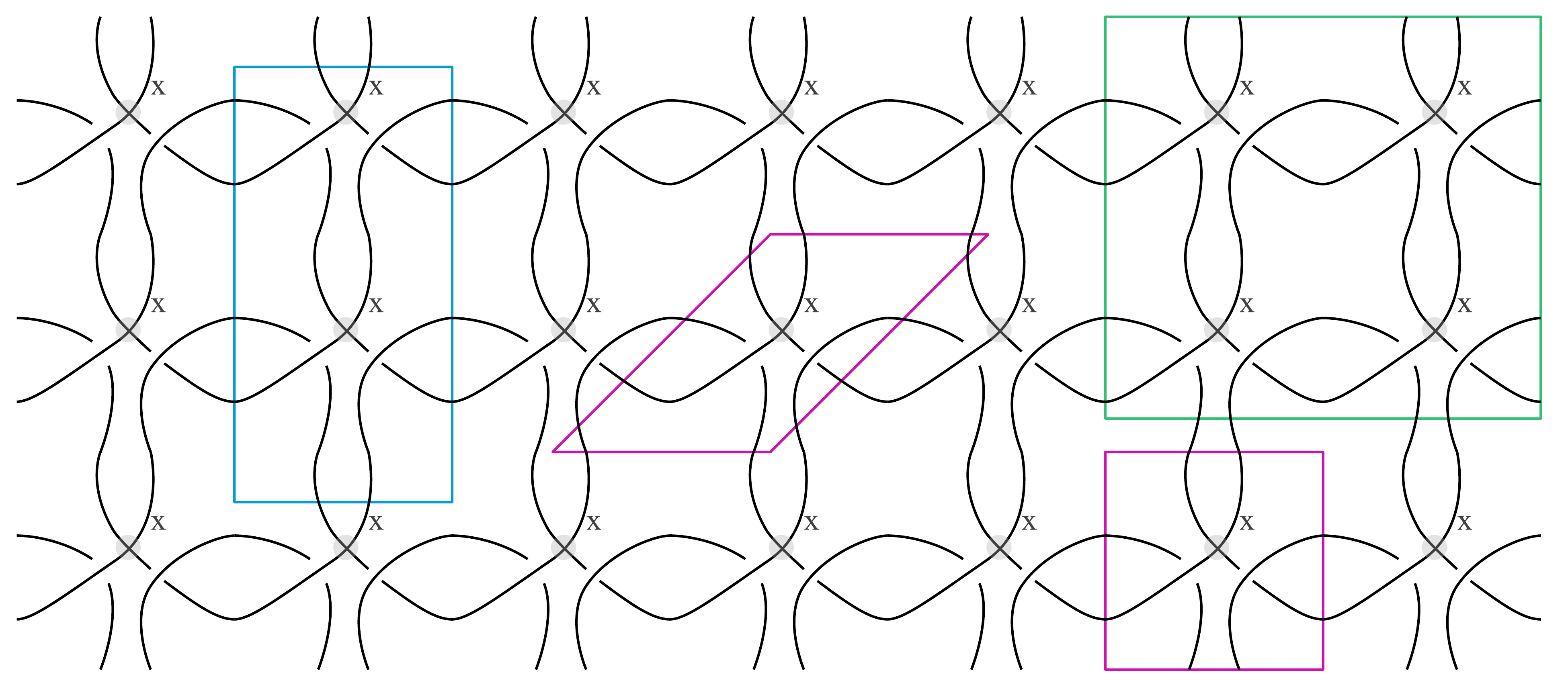}
\caption{\label{minimal-pseudo} Different generating flat pseudo motifs for the same pseudo DP diagram.}
\end{figure}

\section{Pseudo DP tangle equivalence}\label{sec:equivalence}

In this section we study the equivalence of pseudo DP tangles as it is gradually translated on the diagrammatic level of their generating pseudo motifs and flat pseudo motifs, taking into consideration their supporting lattices. 

Our study is in analogy to and uses the equivalence of classical DP tangles, as initiated in \cite{Grishanov1} and established in \cite{DLM1}. A  DP tangle is allowed to undergo \textit{isotopies}, that is, transformations induced by bi-continuous orientation preserving `elastic' deformations (i.e. homeomorphisms) of the thickened plane $\mathbb{E}^2 \times I$. We restrict our scope only to isotopies that preserve the double periodicity. Non-periodic isotopies of a given pseudo DP tangle could be also interesting in the general theory and shall be considered as \textit{defect isotopies}. Two DP tangles are {\it DP isotopic} if they are related via an ambient isotopy of $\mathbb{E}^2 \times I$ that preserves double periodicity. 
This includes \textit{orientation preserving invertible affine transformations} of the plane carrying along the DP tangles. In \cite{DLM1} (see also \cite{Grishanov1}) DP tangle isotopy is translated into combinatorial moves and transformations on the level of generating (flat) motifs and it is referred to as \textit{DP tangle equivalence}. 
We shall expand on the above in the present section in the context of pseudo DP tangles. 

\begin{definition} \label{def:pseudoDPequivalence}
Two pseudo DP tangles are {\it pseudo DP isotopic} if they are related via an ambient isotopy of the thickened plane respecting arc and disc isotopies, that preserves double periodicity. 
\end{definition}

Pseudo DP isotopy extends classical DP isotopy. Recall that any orientation preserving invertible affine transformation of $\mathbb{E}^2$ can be realized as a composition of translations, rotations, shearings and re-scalings. These transformations include both non-area-preserving types, such as re-scalings (e.g. stretchings or contractions), and area-preserving ones, like  translations, rotations, or shear deformations of the plane. Translations, rotations, shearings and re-scaling transformations of the plane, carrying along the DP tangles, give rise on the level of motifs to what we call {\it global equivalence of pseudo  DP tangles}. Arc and disc isotopies in the thickened torus, and torus isotopies give rise to what we call {\it local equivalence of pseudo DP tangles}. On the other hand, a re-scaling of the basis of the plane that affects only the lattice but keeps the DP tangle fixed, translates to a different finite cover or quotient of the (flat) motif. Lattice re-scalings comprise what we call {\it cover equivalence of pseudo  DP tangles}. As we shall see, in the case of pseudo DP tangles there is an important subtlety in the definition of cover equivalence, which led to the necessity of introducing the notion of homologous precrossings. We shall next expand on the above, making the analysis on the motif level, in order to establish pseudo DP tangle equivalence.

\subsection{Local equivalence of (flat) pseudo motifs} \label{sec:Lmotifequivalence}

In the case of classical DP tangles we showed in \cite{DLM1}, using the correspondence between mixed links and toroidal links, that local motif isotopy classes correspond bijectively to mixed link isotopy classes, and  mixed link isotopy is established diagrammatically in \cite[Theorem 5.2]{LR1} (recall Subsection~\ref{sec:mixed}). 
Regarding pseudo DP tangles, {\it local pseudo DP isotopies} are generated by local isotopies of toroidal pseudo links, as described in Subsection~\ref{sec:lift}. On the diagrammatic level, these are induced by surface isotopies of the torus $T^2$ together with the toroidal pseudo Reidemeister equivalence, presented in Subsection~\ref{sec:prel-def}. 

In order to interpret the above local moves on the level of flat pseudo motifs we use the correspondence between mixed pseudo links and toroidal pseudo links discussed in Subsection~\ref{sec:mixed}, in analogy to classical DP tangles. We recall that isotopy for mixed pseudo links is generated by planar isotopy, the classical Reidemeister moves and the pseudo Reidemeister moves for the moving part of an ${\rm H}$-mixed pseudo link diagram $D$  (as exemplified in  Figure~\ref{reid}), together with the extended local isotopies that involve the fixed and the moving part of $D$,  where the fixed components run close and in parallel to the meridian and longitude of the torus $T^2$ respectively. This approach is particularly useful when considering a (flat) pseudo  motif $d$ with respect to the chosen longitude-meridian pair $(l,m)$ of $T^2$, especially when local pseudo isotopy moves take place in $d$ such that some arcs or crossings involved in the moves cross $m$ or $l$ or even both $m$ and $l$. The above lead to the following definition. 
   
\begin{definition} \label{localmoves}
Two flat pseudo motifs are said to be \textit{locally equivalent} if there are related by  {\it local pseudo motif isotopies}, comprising the classical Reidemeister moves, the pseudo Reidemeister moves, and the following surface isotopies: 

\begin{itemize} 
    \item[(a)] planar isotopies within the `interior' of a  flat pseudo motif,
    \smallbreak
    \item[(b)] planar isotopies where an arc before lies within the pseudo motif but the arc afterwards hits one boundary component (the meridian or the longitude), and vice versa,
    \smallbreak
    \item[(c)] planar isotopies where an arc before and the arc afterwards cross one boundary component,
    \smallbreak
    \item[(d)] planar isotopies where an arc before crosses one boundary component but the arc afterwards crosses both boundary components, and vice versa,
    \smallbreak
    \item[(e)] the situation where a classical crossing traverses one boundary component, the meridian or the longitude,
       \smallbreak
    \item[(f)] the situation where a precrossing traverses one boundary component, the meridian or the longitude.
\end{itemize}
\end{definition}

\noindent All local moves are illustrated in  Figure~\ref{fig:iso}. The moves of type f) are the ones needed to complete the pseudo equivalence in addition to the classical case. One of these is illustrated in the bottom right of Figure~\ref{fig:iso}, where the top part shows the move of the pseudo motif, while the bottom part shows the corresponding move of the  flat pseudo motif.

\begin{figure}[ht]
\includegraphics[width=4.5in]{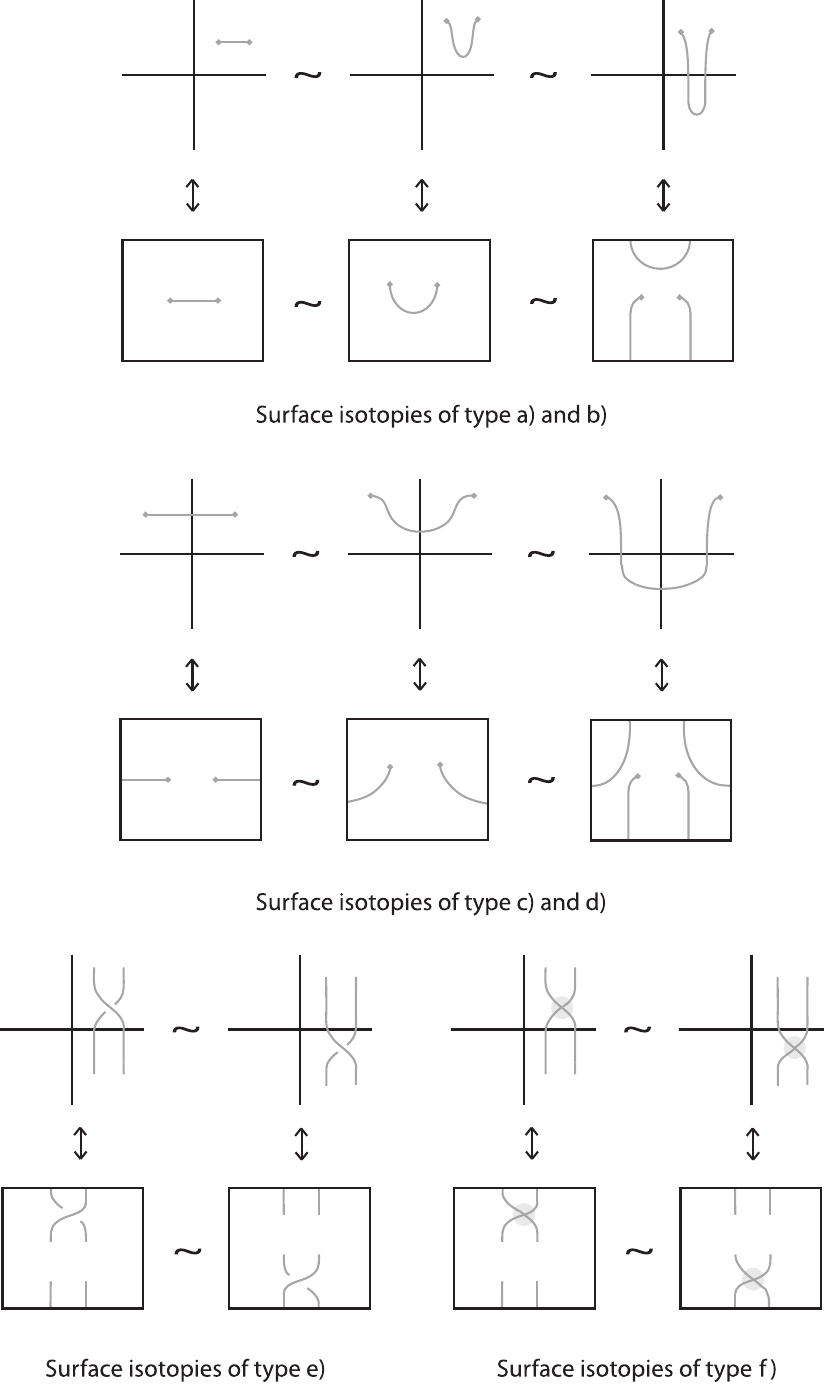}
\caption{\label{fig:iso} Surface isotopy moves.}
\end{figure}

It is worth observing that, if a classical crossing or a precrossing hits the intersection of the longitude-meridian pair, then this move is composite of the ones listed above.

\subsection{Global equivalence of (flat) pseudo motifs}\label{sec:Gmotifequivalence}

Global isotopy  of pseudo DP tangles comprises translations, rotations, re-scalings and shearings of the plane $\mathbb{E}^2$ carrying along the DP tangles. A translation shifts the coordinate system by any vector. Note that a translation of the basis of $\mathbb{E}^2$  leaves the underlying lattice invariant, therefore is not visible on the level of $T^2$. Any other non-integral translation of a pseudo DP tangle  appears on the motif level as a parallel shift of the longitude-meridian pair of $T^2$. Plane rotations preserve the longitude-meridian pair of $T^2$, while rotating the basis of the plane. Lattice re-scalings (stretchings and contractions of the basis vectors) correspond to global  isotopies of $T^2$ (torus inflations and deflations) and they induce re-scalings of the flat motifs.  

Finally, shearings are area-preserving automorphisms of $\mathbb{E}^2$, which on the torus level are related to Dehn twists. A {\it Dehn twist} is an orientation-preserving self-homeomorphism of $T^2$, as described in \cites{Grishanov1,GrishanovP1,GrishanovP2} (for an example see top row of Figure~\ref{Rtwists}). To translate a Dehn twist on the level of pseudo motifs, consider the lattice $\Lambda$ generated by two bases, $B$ and $B'$, as defined in Subsection~\ref{sec:lattice}. We now assign the longitude $l$ of the torus $T^2$ to $u'$ and the meridian $m$ to $v'$ from the basis $B'$. The corresponding covering map, denoted $\rho': \mathbb{E}^2 \rightarrow{} T^2$, creates a new pseudo motif $d'$ that differs from the pseudo motif $d$ associated with the basis $B$ by a finite sequence of Dehn twists. When two flat pseudo motifs are related by a finite sequence of Dehn twists, they are called {\it Dehn equivalent}. For an illustration see bottom row of Figure~\ref{Rtwists}. 

\begin{figure}[H]
\includegraphics[width=5in]{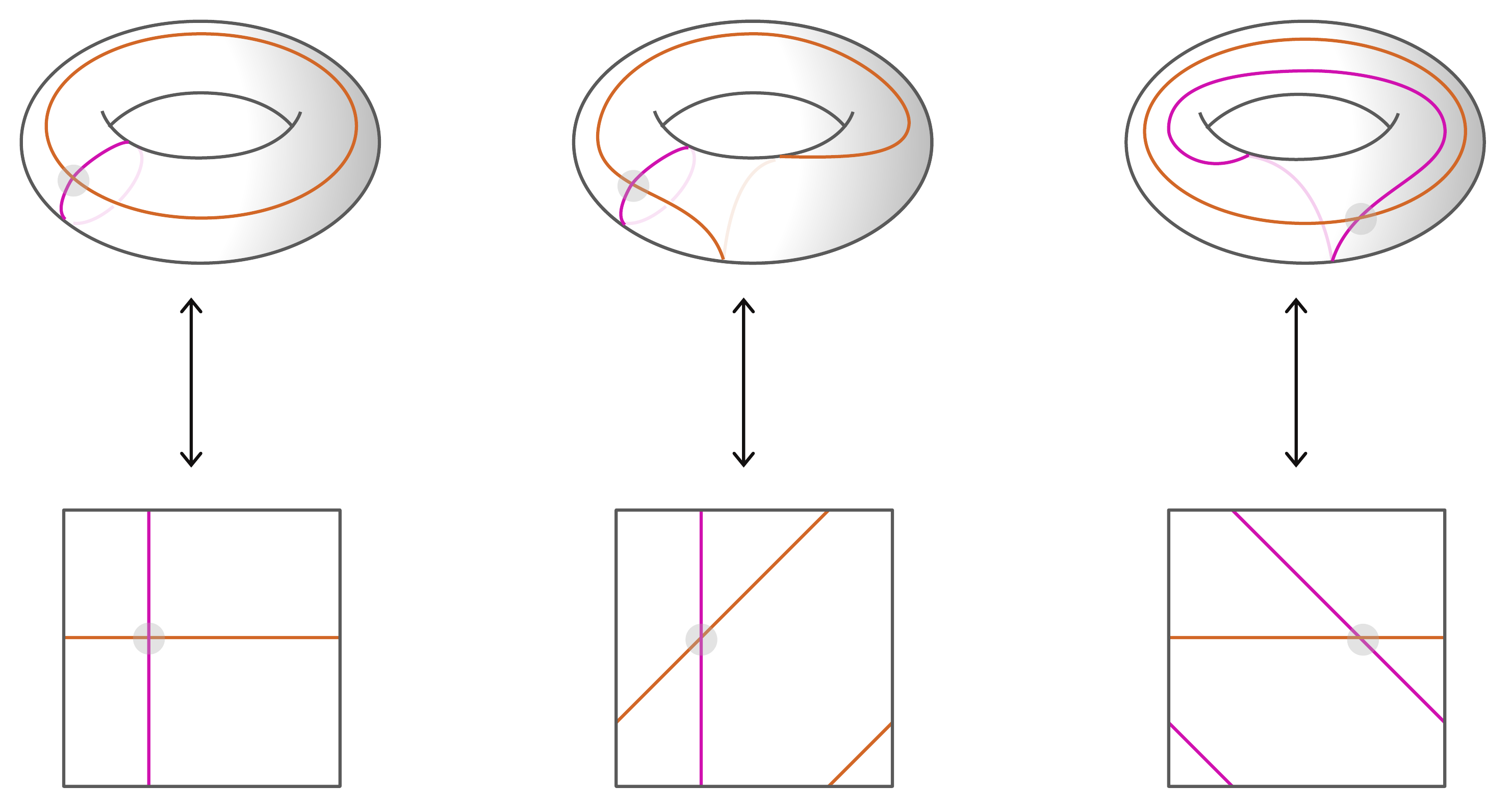}
\caption{\label{Rtwists} 
A pseudo motif (on the left), to which a longitudinal Dehn twist {(center)} and a meridional Dehn twist {(right)} are applied.}
\end{figure}

\begin{definition} \label{globalequiv}
Two flat pseudo motifs are said to be \textit{globally equivalent} if they are related by a finite sequence of  shifts of the longitude-meridian pair of the (flat) torus, rotations, torus inflations and deflations and Dehn twists.
\end{definition}

All global equivalence moves are not affected by the presence of precrossings in the DP tangles, so the theory carries through as in the case of classical DP tangles \cite{DLM1}.

\subsection{Cover equivalence of pseudo DP tangles }\label{sec:scaleequivalence}

Another important move in pseudo DP tangle equivalence is the {\it cover equivalence}, where we have re-scalings of the underlying lattice of a given pseudo DP tangle diagram $d_{\infty}$ keeping the DP tangle fixed. Lattice re-scalings represent distinct finite covers or quotients of the same (flat) pseudo motif. Indeed let $\Lambda_0$, $\Lambda_1$ and $\Lambda_2$ be three periodic lattices (not necessarily distinct) associated with $d_{\infty}$, such that $\Lambda_1 \subseteq \Lambda_0$ and $\Lambda_2 \subseteq \Lambda_0$.  In other words $\Lambda_0$ is a refinement of both $\Lambda_1$ and $\Lambda_2$. Note that the above lattice inclusions are well-defined, since $\Lambda_0$, $\Lambda_1$ and $\Lambda_2$ are given as lattices of the same pseudo DP tangle, so every precrossing belongs to a distinct class of homologous precrossings. Let $d_0 = d_{\infty} / \Lambda_0$, $d_1 = d_{\infty} / \Lambda_1$, and $d_2 = d_{\infty} / \Lambda_2$ represent the three corresponding flat pseudo motifs of $d_{\infty}$. Then, $d_1$ and $d_2$ arise as finite covers of $d_0$. Remark~\ref{rem:homologous} highlights the subtlety of cover equivalence.

\smallbreak
The above lead to the following theorem on pseudo DP tangle equivalence, which is the central result of this paper and which plays a crucial role in defining topological invariants of pseudo DP tangles.

\begin{theorem}[Pseudo DP tangle equivalence]\label{th:equivalence}
Let $\tau_{1,\infty}$ and $\tau_{2,\infty}$ be two pseudo DP tangles in $\mathbb{E}^2 \times I$, with corresponding pseudo DP diagrams $d_{1,\infty}$ and $d_{2,\infty}$. Let also $\Lambda_1$ and $\Lambda_2$ be supporting periodic lattices such that $d_i = d_{i,\infty} / \Lambda_i$ is a flat pseudo motif of $d_{i,\infty}$ for $i \in \{1,2\}$. Then, $\tau_{1,\infty}$ and $\tau_{2,\infty}$ are pseudo DP equivalent if and only if $d_1$ and $d_2$ are related by a finite sequence of local pseudo motif isotopy moves, longitude-meridian shifts, torus inflations and deflations, Dehn twists, and pseudo cover equivalence.
\end{theorem}

\section{Relations with classical DP tangles}\label{sec:resolution}

In this section we present the natural relation between pseudo DP tangles and  classical DP tangles, by means of the resolution sets. Building on the foundational ideas introduced for planar and toroidal pseudo links in \cite{DLM4} (recall Section~\ref{sec:were}), we adapt the resolution and the weighted resolution sets  to the framework of  pseudo DP tangles. This extension allows one to study  pseudo DP tangles via classical  DP tangles and invariants of theirs.

\begin{definition}\label{def:ReSet}
A {\it resolution} of a  pseudo DP diagram $d_\infty$ is a classical DP tangle obtained by assigning to every precrossing belonging to the same class of homologous precrossings, the same  type of crossing  (positive or negative). The set of all resolutions of $d_\infty$  is its {\it resolution set}, denoted $Re(d_\infty)$. Further, let $d$ be a (flat) pseudo motif diagram of  $d_{\infty}$. A {\it resolution} of $d$ and the {\it resolution set} of $d$ are defined as in Definition~\ref{wrstormixed} for toroidal pseudo link diagrams.
\end{definition} 

It follows immediately from Definition~\ref{def:ReSet} that a resolution of $d_\infty$ arises equivalently as the universal cover of the respective resolution of $d$. This is also the case for any finite cover $\tilde{d}$ of $d$, since it retains the same (hidden) information on homologous precrossings, which resolve in the same crossing type. Therefore, every element of the resolution set of the finite cover $\tilde{d}$ of $d$ is a finite cover of an element of the resolution set of $d$. Figure~\ref{tortref} top illustrates an example of the resolution set of a flat pseudo motif diagram containing a precrossing $x$. Figure~\ref{tortref} bottom illustrates a double cover of the flat pseudo motif diagram of the top and its resolution set. Hence, the pseudo DP tangle generated by either of these flat pseudo motifs has a resolution set of two elements.

\begin{figure}[ht]
\begin{center}
\includegraphics[width=3.8in]{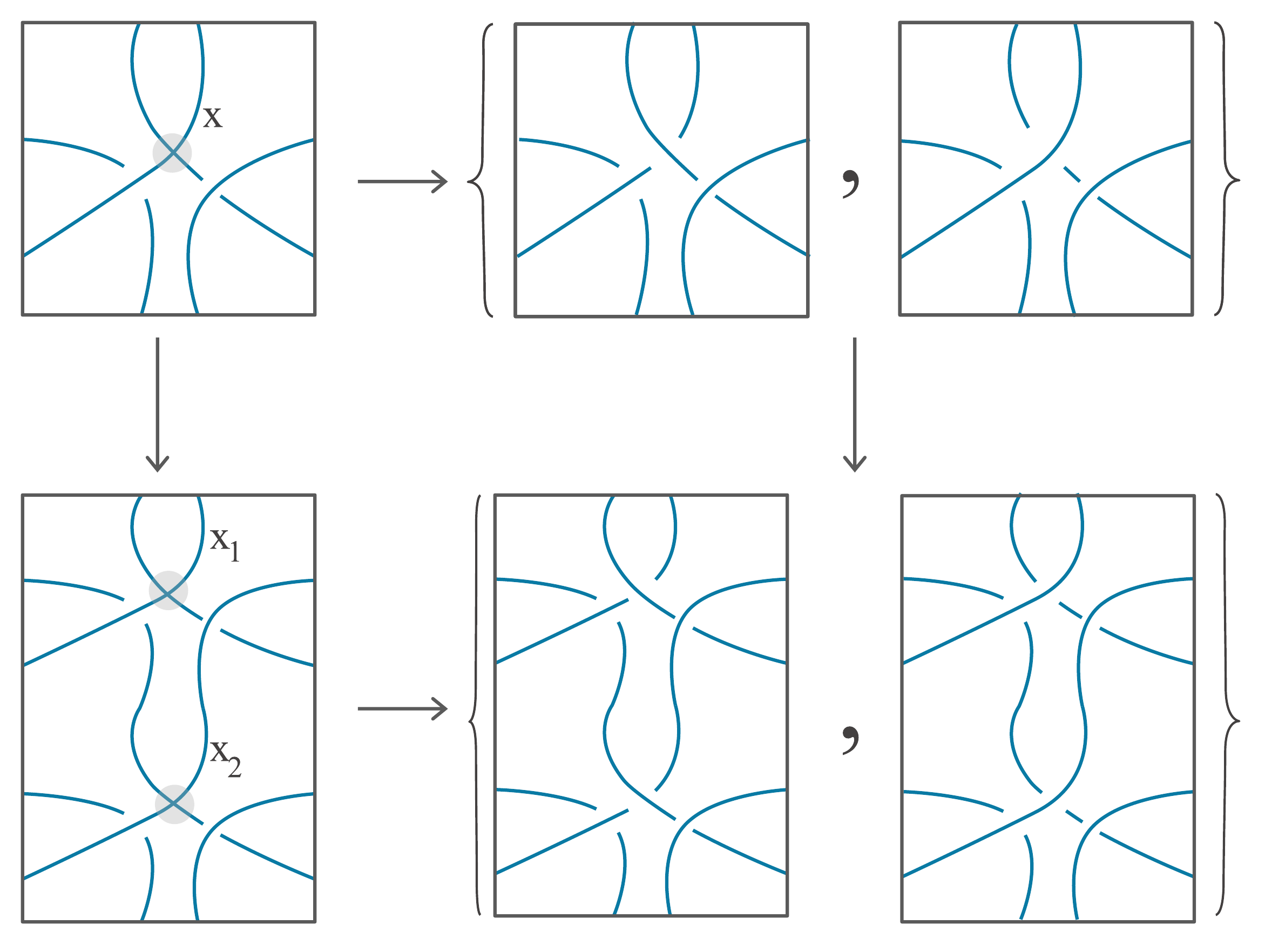}
\end{center}
\caption{Top: the resolution set of a  flat pseudo motif diagram. Bottom: their respective double covers.}
\label{tortref}
\end{figure}

Note, however, that if a (flat) pseudo motif $\tilde{d}'$ is identical to the finite cover  $\tilde{d}$  of $d$, except that it has no marking of homologous precrossings (see Remark~\ref{rem:homologous}), then each of the precrossings of $\tilde{d}'$ must resolve independently. Figure~\ref{tortref2} illustrates an example of the subtlety of cover equivalence discussed above, where the precrossings $x$ and $y$ are not marked as homologous. So, the resolution set of the motif $\tilde{d}'$ contains four elements and thus differs from the  resolution set of $\tilde{d}$ in Figure~\ref{tortref} which contains two elements. 

\begin{figure}[H]
\begin{center}
\includegraphics[width=5.05in]{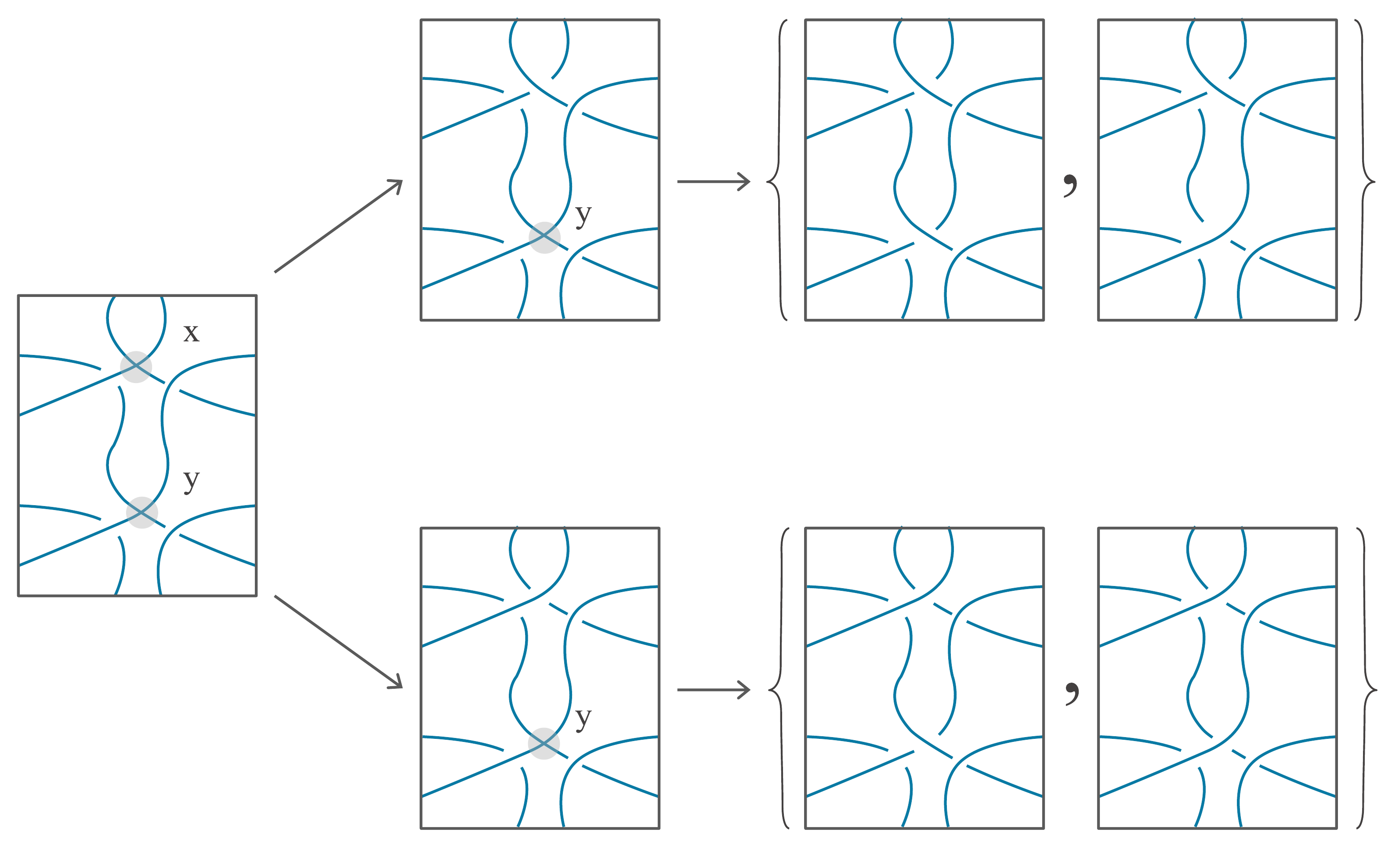}
\end{center}
\caption{The resolution set of a flat pseudo motif diagram with two non homologous precrossings.}
\label{tortref2}
\end{figure}

By the theory of DP tangle equivalence (recall Introduction and \cite{DLM1}) and  that of pseudo DP tangles (Theorem~\ref{th:equivalence}) it is straightforward to verify that the resolution set of a pseudo DP diagram $d_\infty$  remains invariant under pseudo DP tangle  equivalence, and hence,  it is an invariant of $d_\infty$. 
 
\begin{theorem}\label{th:ReDP}
The resolution set $Re(d_\infty)$ is a topological invariant of a pseudo DP diagram $d_\infty$. Furthermore, applying any invariant of classical DP tangles to the elements of  $Re(d_\infty)$ also induces an invariant for $d_\infty$. 
\end{theorem}

Finally, in analogy to the theory of toroidal pseudo links (recall Section~\ref{sec:were}) we define: 

\begin{definition}\label{def:WeReSet}
The {\it weighted resolution set} or {\it WeRe set} of a  pseudo DP diagram $d_\infty$ is a collection of ordered pairs \((d^i_\infty, p_{i})\), where $d^i_\infty$ represents a resolution of $d_\infty$, and \( p_{i} \) denotes the probability of obtaining from $d_\infty$ the equivalence class of $d^i_\infty$  through a random assignment of crossing types to homologous precrossings, with equal likelihood for positive and negative crossings. 
\end{definition}

We conclude the paper with the following theorem which follows immediately from Definitions~\ref{def:WeReSet} and~\ref{def:ReSet} and Theorem~\ref{th:ReDP}:

\begin{theorem}\label{th:WeReDP}
The WeRe set is a topological invariant of pseudo DP tangles.
\end{theorem}


\begin{thebibliography}{99}

\bibitem{D} I. Diamantis, Pseudo links and singular links in the Solid Torus, {\it Communications in Mathematics}, Vol. 31, Issue 1 (2023).

\bibitem{DLM1} I. Diamantis, S. Lambropoulou, S. Mahmoudi, Equivalences of doubly periodic tangles, \textit{arXiv:2310.00822} (2023).

\bibitem{DLM4} I. Diamantis, S. Lambropoulou, S. Mahmoudi, From annular to toroidal pseudo knots, {\em Symmetry} {\bf 2024}, {\it 16}(10), 1360. https://doi.org/10.3390/sym16101360.

\bibitem{H} R. Hanaki, Pseudo diagrams of knots, links and spatial graphs, {\it Osaka J. Math.}, {\bf 47}, (2010) 863--883.

\bibitem{HJMR} A. Henrich, R. Hoberg, S. Jablan, L. Johnson, E. Minten, L. Radovic, The theory of pseudoknots, {\it J. Knot Theory Ramif. }, {\bf 22}, No. 07, (2013) 1350032.

\bibitem{Grishanov1} S. Grishanov, V. Meshkov, A. Omelchenko, Kauffman-type polynomial invariants for doubly periodic structures, {\it J. Knot Theory Ramif.} \textbf{16} (2007), 779--788. 

\bibitem{GrishanovP1} S. Grishanov, V. Meshkov, A. Omelchenko, A topological study of textile structures. Part I: an introduction to topological methods. {\it Text. Res. J.}. {\bf 79} (2009),702--713.

\bibitem{GrishanovP2} S. Grishanov, V. Meshkov, A. Omelchenko, A topological study of textile structures. Part II: topological invariants in application to textile structures, {\it Text. Res. J.} \textbf{79} (2009), 822--836.

\bibitem{LR1} S. Lambropoulou, C.P. Rourke, Markov's theorem in $3$-manifolds, \emph{Topology and its Applications} {\bf 78},
(1997) 95-122.

\bibitem{Yaghi} Y. Liu, M. O'Keeffe, M.M.J. Treacy, O.M. Yaghi, The geometry of periodic knots, polycatenanes and weaving from a chemical perspective: A library for reticular chemistry. \textit{Chem. Soc. Rev.} {\bf 2018}, {\it 47}, 4642--4664.
	
\bibitem{Treacy} M. O'Keeffe, M.M.J. Treacy, Crystallographic descriptions of regular 2-periodic weavings of threads, loops and nets. \textit{Acta Cryst. A} {\bf 2020},  {\it 76}, 110--120.

\bibitem{Sonia} S. Mahmoudi, On the classification of periodic weaves and universal cover of links in thickened surfaces, {\em Commun. Korean Math. Soc.}, {\bf 39} (2024), No. 4, pp. 997-1025.

\end{thebibliography}
\end{document}